\documentclass[11pt]{article}
\usepackage[T1]{fontenc}
\usepackage{lmodern,amsmath,amssymb,amsthm,mathtools,microtype}
\usepackage[margin=1in]{geometry}
\usepackage[colorlinks=true,linkcolor=blue,citecolor=blue,urlcolor=blue]{hyperref}
\hypersetup{pdftitle={An interpretation of the BMO boundedness of Hormander-Mikhlin multipliers}}
\newtheorem{theorem}{Theorem}[section]
\newtheorem{lemma}[theorem]{Lemma}
\newtheorem{corollary}[theorem]{Corollary}
\theoremstyle{remark}\newtheorem{remark}[theorem]{Remark}
\newcommand{\RR}{\mathbb R}
\newcommand{\HH}{\mathcal H}
\newcommand{\cM}{\mathcal M}
\newcommand{\bmo}{\mathrm{bmo}}
\newcommand{\BMO}{\mathrm{BMO}}
\newcommand{\norm}[1]{\left\|#1\right\|}
\newcommand{\esssup}{\operatorname*{ess\,sup}}
\title{An interpretation of the BMO boundedness of H\"ormander--Mikhlin multipliers}
\author{Tao Mei}
\date{}
\begin{document}
\maketitle
\vspace{-2em}
\begin{abstract}
We prove a pointwise domination of Poisson variances for H\"ormander--Mikhlin
multipliers on the real line. The estimate holds in operator order for
operator-valued functions and gives a direct Fourier-side proof of BMO
boundedness. The same argument gives an algebraic proof of the
one-dimensional Schur multiplier estimate
of Conde-Alonso, Gonz\'alez-P\'erez, Parcet and Tablate
(\emph{Annals of Mathematics}, 2023), formulated here using Poisson BMO.
The endpoint estimates are proved for $L^2$ functions and Hilbert--Schmidt
operators; interpolation gives the corresponding $L^p$ and Schatten estimates.
\end{abstract}

\section{Introduction}\label{sec:introduction}
The H\"ormander--Mikhlin theorem provides a sufficient condition for the
$L^p$ boundedness of Fourier multipliers. In one dimension, it says that
the multiplier $T_\mu$ with symbol $\mu$ satisfying
\[
 M=\|\mu\|_\infty<\infty,\qquad
 D=\esssup_{\xi\ne0}|\xi\mu'(\xi)|<\infty,
\]
with local absolute continuity on each open half-line, is bounded on
$L^p$ for all $1<p<\infty$. These multipliers are usually treated within
Calder\'on--Zygmund theory. One standard route combines an
$L^\infty\to\mathrm{BMO}$ endpoint estimate with $L^2$ boundedness,
interpolation, and duality.

This article gives a direct Fourier-side proof of the Poisson BMO
estimate without estimating the singular-integral kernel of the multiplier.
We use the variance definition of Poisson semigroup BMO:
\[
 \|f\|_{\bmo(P)}
 =\sup_{t>0}\|P_t|f|^2-|P_tf|^2\|_\infty^{1/2},
 \qquad \widehat{P_tf}(\xi)=e^{-2\pi t|\xi|}\widehat f(\xi).
\]
For scalar functions on $\RR$, it is equivalent to the quadratic interval
$\BMO_2$ seminorm; a self-contained proof is given in
Appendix~\ref{app:interval-bmo}. For matrix-valued and operator-valued
functions we use both column and row variances.
It is well known that the Hilbert transform preserves Poisson variance.
Our Fourier-side argument gives the following pointwise comparison for
H\"ormander--Mikhlin multipliers:
\[
 P_t|T_\mu f|^2-|P_tT_\mu f|^2
 \le(M+2D)^2\bigl(P_t|f|^2-|P_tf|^2\bigr),
\]
for $f\in L^2(\RR)$ and the Fourier multiplier $T_\mu$ with symbol $\mu$.
The same proof works in operator order for operator-valued $L^2$ functions.

For Schur multipliers on $S_2(L^2(\RR))$, the Poisson semigroup multiplies a kernel
by $e^{-t|x-y|}$. Our argument provides an algebraic proof of
the one-dimensional version of the Schur multiplier endpoint estimate of
Conde-Alonso, Gonz\'alez-P\'erez, Parcet and Tablate
\cite[Section 2 and Remark 1.6]{CGPT}, formulated here for Poisson BMO.
 
Section~\ref{sec:fourier} proves the scalar Fourier theorem and its
operator-valued extension.
Section~\ref{sec:schur} proves the Schur multiplier theorem directly on
kernels. The BMO endpoint proofs are self-contained.  Appendix~\ref{app:interval-bmo} records
the standard equivalence between Poisson semigroup BMO and quadratic
interval $\BMO_2$, including its matrix-valued version.
Appendix~\ref{app:interpolation} applies BMO--$L^2$ interpolation and its
Schatten counterpart. Appendix~\ref{app:hardy} gives the direct Hardy-space estimate. For the passage to Schatten classes, we use the
semigroup BMO framework of Mei \cite{MeiTent} and the interpolation theory
of Junge--Mei \cite{JM}.

\section{H\"ormander--Mikhlin theory on the real line}\label{sec:fourier}

Let $P_t f=p_t*f$ be the Poisson semigroup on $\RR$, with
$p_t(x)=t/(\pi(t^2+x^2))$. We use the Fourier convention
\[
 \widehat f(\xi)=\int_\RR f(x)e^{-2\pi i x\xi}\,dx,\qquad
 f(x)=\int_\RR\widehat f(\xi)e^{2\pi i x\xi}\,d\xi
\]
for Schwartz functions and its Plancherel extension on $L^2(\RR)$.
Thus $\widehat{P_tf}(\xi)=e^{-2\pi t|\xi|}\widehat f(\xi)$.
For $f\in L^2(\RR,(1+x^2)^{-1}\,dx)$ put
\[
 V_t^P(f)(x)=P_t|f|^2(x)-|P_tf(x)|^2,\qquad
 \|f\|_{\bmo(P)}=\sup_{t>0}\|V_t^P(f)\|_\infty^{1/2}.
\]
The Poisson averages are finite and continuous for each $t>0$, since
$p_t(x-y)$ is comparable to $(1+y^2)^{-1}$ for fixed $t,x$, locally
uniformly in $x$. The space $\bmo(P)$ consists of these weighted
$L^2$ functions with finite seminorm, modulo constants.
We first establish the Fourier identities and estimates for
$f\in L^2(\RR)$, allowing its BMO seminorm to be infinite.
We begin with scalar-valued functions and then treat matrix-valued and
operator-valued functions.

\subsection{Covariance factorization and Hardy inequality}
For $u,v>0$, direct integration gives
\begin{equation}\label{eq:continuous-scalar}
 e^{-t|u-v|}-e^{-t(u+v)}
 =2t\int_0^{\min(u,v)}e^{-t(u-s)}e^{-t(v-s)}\,ds.
\end{equation}
For a Hilbert-space-valued $h\in L^2(\RR_+)$, the substitution $u=sv$
and Minkowski's inequality give the Hardy estimate
\begin{equation}\label{eq:hardy}
 \begin{aligned}
 \left\|s\mapsto\int_s^\infty\frac{h(u)}u\,du\right\|_{L^2(ds)}
 &\le\int_1^\infty\|h(sv)\|_{L^2(ds)}\frac{dv}{v}\\
 &=\|h\|_2\int_1^\infty v^{-3/2}\,dv
 =2\|h\|_2.
 \end{aligned}
\end{equation}
In particular, an integral kernel bounded in absolute value by
$D\mathbf1_{\{u>s\}}/u$ has $L^2$ operator norm at most $2D$,
also for Hilbert-space-valued inputs.

The following lemma expresses each Poisson variance as an $L^2$
integral of squared Fourier tails. This factorization is the key ingredient
of the argument.

\begin{lemma}[Fourier-tail identity]\label{lem:fourier-tails}
Let $f\in L^2(\RR)$, with Fourier transform defined by Plancherel.
For $t>0$, $x\in\RR$, $\sigma=\pm1$ and $s>0$, define
\[
 F_s^\sigma(t,x)=\int_s^\infty\widehat f(\sigma u)
 e^{-2\pi t(u-s)}e^{2\pi i\sigma ux}\,du.
\]
Then
\begin{equation}\label{eq:real-variance}
 V_t^P(f)(x)=4\pi t\sum_{\sigma=\pm1}\int_0^\infty
 |F_s^\sigma(t,x)|^2\,ds.
\end{equation}
\end{lemma}
\begin{proof}
First suppose that $f\in\mathcal S(\RR)$. Then $\widehat f\in L^1(\RR)$,
and Fubini's theorem gives
\begin{align*}
 P_t|f|^2(x)
 &=\int_\RR\int_\RR
 \overline{\widehat f(\xi)}\widehat f(\eta)
 e^{2\pi i x(\eta-\xi)}e^{-2\pi t|\eta-\xi|}\,d\xi\,d\eta,\\
 |P_tf(x)|^2
 &=\int_\RR\int_\RR
 \overline{\widehat f(\xi)}\widehat f(\eta)
 e^{2\pi i x(\eta-\xi)}e^{-2\pi t(|\xi|+|\eta|)}\,d\xi\,d\eta.
\end{align*}
Subtracting these expressions yields the Fourier representation of
$V_t^P(f)(x)$. When $\xi$ and $\eta$ have opposite signs,
$|\eta-\xi|=|\xi|+|\eta|$, so their contribution vanishes.
For the remaining terms, write $\xi=\sigma u$ and $\eta=\sigma v$,
where $u,v>0$ and $\sigma=\pm1$. Thus
\[
 V_t^P(f)(x)=\sum_{\sigma=\pm1}\int_0^\infty\int_0^\infty
 \overline{\widehat f(\sigma u)}\widehat f(\sigma v)
 e^{2\pi i\sigma x(v-u)}\bigl(e^{-2\pi t|u-v|}-e^{-2\pi t(u+v)}\bigr)\,du\,dv.
\]
Apply \eqref{eq:continuous-scalar} with $t$ replaced by $2\pi t$. The region
$0<s<\min(u,v)$ is equivalently described by $s>0$, $u>s$ and $v>s$.
Changing the order of integration therefore gives
\begin{align*}
 V_t^P(f)(x)
 &=4\pi t\sum_{\sigma=\pm1}\int_0^\infty
 \left(\int_s^\infty\overline{\widehat f(\sigma u)}
 e^{-2\pi t(u-s)}e^{-2\pi i\sigma ux}\,du\right)\\
 &\hspace{38mm}\times
 \left(\int_s^\infty\widehat f(\sigma v)
 e^{-2\pi t(v-s)}e^{2\pi i\sigma vx}\,dv\right)\,ds\\
 &=4\pi t\sum_{\sigma=\pm1}\int_0^\infty
 \overline{F_s^\sigma(t,x)}F_s^\sigma(t,x)\,ds,
\end{align*}
which is \eqref{eq:real-variance}.

To justify the last interchange, take absolute values before integrating
in $s$. By the same factorization, the resulting triple integral is
\[
 \sum_{\sigma=\pm1}\int_0^\infty\int_0^\infty
 |\widehat f(\sigma u)|\,|\widehat f(\sigma v)|
 \bigl(e^{-2\pi t|u-v|}-e^{-2\pi t(u+v)}\bigr)\,du\,dv\le\|\widehat f\|_1^2,
\]
because $0\le e^{-2\pi t|u-v|}-e^{-2\pi t(u+v)}\le1$. Hence all the preceding interchanges
are absolutely justified; in particular, the Fourier tails are square
integrable in $s$.

For general $f\in L^2(\RR)$, each tail integral converges absolutely
by Cauchy--Schwarz. The substitution $u=s+r$ and Minkowski's
inequality give
\[
 \|F_\cdot^\sigma(t,x)\|_{L^2(ds)}
 \le\int_0^\infty e^{-2\pi tr}
       \|\widehat f(\sigma(s+r))\|_{L^2(ds)}\,dr
 \le(2\pi t)^{-1}\|\widehat f(\sigma\cdot)\|_{L^2(\RR_+)}.
\]
Choose $f_n\in\mathcal S(\RR)$ converging to $f$ in $L^2$.
The displayed bound and Plancherel show that their tails converge
in $L^2(ds)$ to those of $f$, uniformly in $x$. Moreover,
\begin{align*}
 \|P_t(|f_n|^2-|f|^2)\|_\infty
 &\le\|p_t\|_\infty\||f_n|^2-|f|^2\|_1\longrightarrow0,\\
 \|P_t(f_n-f)\|_\infty
 &\le\|p_t\|_2\|f_n-f\|_2\longrightarrow0.
\end{align*}
Thus both sides of \eqref{eq:real-variance} converge, proving the
identity for every $t>0$ and $x\in\RR$.
\end{proof}

\subsection{The scalar theorem}
\begin{theorem}\label{thm:real-fourier}
Let $\mu\in L^\infty(\RR)$ be locally absolutely continuous on each
open half-line, and set
\[
 M=\|\mu\|_\infty,\qquad
 D=\esssup_{\xi\ne0}|\xi\mu'(\xi)|<\infty.
\]
The Fourier multiplier $\widehat{T_\mu f}=\mu\widehat f$ satisfies, for
every $f\in L^2(\RR)$,
\[
 V_t^P(T_\mu f)(x)\le(M+2D)^2 V_t^P(f)(x),\qquad t>0,\ x\in\RR.
\]
Consequently,
\[
\|T_\mu f\|_{\bmo(P)}\le(M+2D)\|f\|_{\bmo(P)},
 \qquad f\in L^2(\RR)\cap\bmo(P).
\]
In particular,
\[
 \|T_\mu f\|_{\bmo(P)}\le(M+2D)\|f\|_\infty,
 \qquad f\in L^2(\RR)\cap L^\infty(\RR).
\]
\end{theorem}

\begin{proof}
Fix $t>0$ and $x\in\RR$. Let $F_s^\sigma$ be the tails from
Lemma~\ref{lem:fourier-tails}, and let $G_s^\sigma$ be the corresponding
tails of $T_\mu f$. Both $f$ and $T_\mu f$ belong to
$L^2(\RR)$, so the lemma applies to both functions.
The tails satisfy
\[
 G_s^\sigma=\mu(\sigma s)F_s^\sigma
 +\int_s^\infty\sigma\mu'(\sigma u)e^{-2\pi t(u-s)}F_u^\sigma\,du.
\]
Indeed, insert
$\mu(\sigma v)=\mu(\sigma s)+\int_s^v\sigma\mu'(\sigma u)\,du$
into the defining integral and use Fubini. The diagonal term has norm
at most $M$ on the tail space, and \eqref{eq:hardy} bounds the integral
term by $2D$. For completeness, the Fubini interchange for each $s>0$
is justified by
\[
 D\int_s^\infty |\widehat f(\sigma v)|e^{-2\pi t(v-s)}\log(v/s)\,dv<\infty,
\]
which follows from Cauchy--Schwarz. Hence
\[
 \sum_\sigma\int_0^\infty|G_s^\sigma|^2\,ds
 \le(M+2D)^2\sum_\sigma\int_0^\infty|F_s^\sigma|^2\,ds.
\]
Identity \eqref{eq:real-variance} proves the variance estimate. Taking
suprema gives the BMO bound, and
$\|f\|_{\bmo(P)}\le\|f\|_\infty$ gives the last assertion.
\end{proof}

\subsection{Matrix-valued and operator-valued functions}
Let $(\cM,\tau)$ be a semifinite von Neumann algebra, and write
\[
 L^2(\RR;\cM)=L^2(\RR;L^2(\cM,\tau)).
\]
This is a Hilbert space, so the Fourier transform and $T_\mu$ act by
Plancherel. For $\cM=M_q$ with its usual trace, it is the space of
matrix-valued functions with square-integrable entries. For
$\cM=B(\HH)$ with the usual trace, the values are Hilbert--Schmidt operators.

For $f\in L^2(\RR;\cM)$ set
\begin{align*}
 V_{t,c}^P(f)(x)&=P_t(f^*f)(x)-(P_tf(x))^*P_tf(x),\\
 \|f\|_{\bmo_c(P)}&=\sup_{t>0}\|V_{t,c}^P(f)\|_\infty^{1/2},\\
 \|f\|_{\bmo_r(P)}&=\|f^*\|_{\bmo_c(P)},\qquad
 \|f\|_{\bmo(P)}=\max\{\|f\|_{\bmo_c(P)},\|f\|_{\bmo_r(P)}\}.
\end{align*}
For each $t,x$, the variance is a positive element of $L^1(\cM,\tau)$.
Its operator norm may be infinite. Here $P_t$ acts only on the real
variable, and the essential supremum uses the operator norm of $\cM$.

\begin{corollary}\label{cor:matrix-fourier}
Under the hypotheses of Theorem~\ref{thm:real-fourier}, every
$f\in L^2(\RR;\cM)$ satisfies
\[
 V_{t,c}^P(T_\mu f)(x)\le(M+2D)^2 V_{t,c}^P(f)(x),
 \qquad t>0,\ x\in\RR,
\]
in positive $L^1(\cM,\tau)$ order, and the same assertion holds for
the row variance. Consequently,
\[
 \|T_\mu f\|_{\bmo(P)}\le(M+2D)\|f\|_{\bmo(P)},
 \qquad f\in L^2(\RR;\cM)\cap\bmo(P).
\]
For $f\in L^2\cap L^\infty$, the right-hand side is at most
$(M+2D)\|f\|_\infty$. The constant is independent of $\cM$ and,
in particular, of the matrix size.
\end{corollary}
\begin{proof}
Initially let $f=\sum_{j=1}^N f_ja_j$, where the $f_j$ are scalar
Schwartz functions and $a_j\in\cM\cap L^2(\cM,\tau)$.
Expanding in the same order as in Lemma~\ref{lem:fourier-tails} gives
\[
 V_{t,c}^P(f)(x)
 =4\pi t\sum_\sigma\int_0^\infty(F_s^\sigma)^*F_s^\sigma\,ds.
\]
No coefficient products have been interchanged. Apply the Hilbert-space
Hardy inequality \eqref{eq:hardy} to $F_s^\sigma\eta$, for each vector $\eta$ in a
concrete representation of $\cM$. This proves column variance domination
in operator order. For the row estimate use
$(T_\mu f)^*=T_{\mu^\sharp}(f^*)$, where
$\mu^\sharp(\xi)=\overline{\mu(-\xi)}$ has the same constants.

For general $f$, approximate by such finite sums $f_n$ in $L^2(\RR;\cM)$.
For each $t,x$,
\[
 \|P_t(f_n^*f_n-f^*f)(x)\|_{L^1(\cM)}
 \le\|p_t\|_\infty\|f_n-f\|_2(\|f_n\|_2+\|f\|_2)
 \longrightarrow0.
\]
Also $P_tf_n(x)\to P_tf(x)$ in $L^2(\cM,\tau)$, hence their squared
moduli converge in $L^1(\cM,\tau)$. The same applies to $T_\mu f_n$,
by its $L^2$ boundedness. The positive cone of $L^1$ is closed, so both
variance inequalities pass to the limit. The argument commutes with
matrix amplification.

The calculation on finite sums uses only bounded operator coefficients;
it therefore also applies to such sums in any concrete von Neumann algebra.
\end{proof}

\begin{remark}
Heat-semigroup and higher-dimensional extensions will appear in \cite{CMR}.
\end{remark}

\begin{remark}[Duality and weak-* density]\label{rem:hardy-duality}
The boundedness of $T_\mu$ on the classical Poisson Hardy space
$H_P^1(\RR)$ already follows from the BMO boundedness by
$H^1$--BMO duality and the weak-* density of Schwartz functions in BMO
modulo constants. The bounded form of this density is useful here:
every $b\in\bmo(P)$ is a weak-* limit of $b_j\in\mathcal S(\RR)$ with
$\sup_j\|b_j\|_{\bmo(P)}\le C\|b\|_{\bmo(P)}$.
Indeed, let $\eta_R(x)=\eta(x/R)$, where $\eta\in C_c^\infty$ is $1$
on $[-1,1]$ and supported in $[-2,2]$, and let $\varphi_\varepsilon$
be a smooth compactly supported approximate identity. The functions
\[
 b_{R,\varepsilon}
 =\varphi_\varepsilon*
       \bigl[\eta_R(b-b_{(-R,R)})\bigr]
 \in\mathcal S(\RR)
\]
have uniformly controlled quadratic interval BMO norms. This follows from
$|\eta_R'|\lesssim R^{-1}$ and the logarithmic bound for differences of
interval averages, obtained by telescoping as in
Appendix~\ref{app:interval-bmo}. Mollification preserves the bound.
Pairing first against compactly supported smooth mean-zero functions shows
convergence as $R\to\infty$, $\varepsilon\to0$; their density in $H^1$
and the uniform BMO bound give weak-* convergence. The equivalence in
Appendix~\ref{app:interval-bmo} transfers this assertion to $\bmo(P)$.

Consequently, duality \cite{FS} allows the Hardy norm to be tested on
Schwartz functions:
\[
 \|h\|_{H_P^1}\lesssim
 \sup_{\substack{b\in\mathcal S(\RR)\\\|b\|_{\bmo(P)}\le1}}
     |\langle h,b\rangle|,
 \qquad h\in H_P^1.
\]
For $f$ in the dense core
$\mathcal S_0=\{f:\widehat f\in C_c^\infty(\RR\setminus\{0\})\}$,
Plancherel gives
\[
 \langle T_\mu f,b\rangle
 =\langle f,T_{\overline\mu}b\rangle.
\]
Theorem~\ref{thm:real-fourier}, applied to $\overline\mu$, therefore yields
\[
 |\langle T_\mu f,b\rangle|
 \le C(M+2D)\|f\|_{H_P^1}\|b\|_{\bmo(P)}.
\]
Here $T_\mu f\in H_P^1$ on this core: its Fourier transform and that of
its Hilbert transform are compactly supported $W^{1,2}$ functions,
so weighted Cauchy--Schwarz makes both inverse transforms integrable.
Taking the supremum and extending by density proves the Hardy bound.
Applying this bound to $\overline\mu$ and taking the adjoint gives
\[
 \widetilde T_\mu
 =\bigl(T_{\overline\mu}|_{H_P^1}\bigr)^*
 :\bmo(P)\longrightarrow\bmo(P),
 \qquad
 \|\widetilde T_\mu\|\le C(M+2D).
\]
This is the unique bounded weak-* continuous extension of the original
multiplier. It agrees with $T_\mu$ on $L^2\cap\bmo(P)$ by Plancherel
pairing against $\mathcal S_0$ and density in $H_P^1$.
We include a direct Fourier-side proof of the $H_P^1$ boundedness in
Appendix~\ref{app:hardy}, using the standard atomic characterization in
place of duality.
\end{remark}

\section{Schur multipliers}\label{sec:schur}

Let $H=L^2(\RR)$ and identify $S_2(H)$ with square-integrable kernels. The Schur--Poisson semigroup is
\[
 K_{S_tX}(x,y)=e^{-t|x-y|}K_X(x,y).
\]
It extends to a normal unital completely positive semigroup on $B(H)$: if $(U_uf)(x)=e^{iux}f(x)$ and $p_t(u)=t/(\pi(t^2+u^2))$, then
\[
 S_tX=\int_\RR p_t(u)U_uXU_u^*\,du.
\]
For $X\in S_2(H)$ define
\[
 \norm{X}_{\bmo_c(S)}=\sup_{t>0}\norm{S_t(X^*X)-(S_tX)^*S_tX}^{1/2},
\]
and set $\|X\|_{\bmo_r(S)}=\|X^*\|_{\bmo_c(S)}$ and
$\|X\|_{\bmo(S)}=\max\{\|X\|_{\bmo_c(S)},\|X\|_{\bmo_r(S)}\}$.
These are Poisson BMO seminorms; their null spaces are factored out when
taking completions. Complete positivity gives
$\|X\|_{\bmo(S)}\le\|X\|$. A bounded measurable scalar symbol $m$ defines $M_m$ on $S_2(H)$ by $K_{M_mX}=mK_X$.

\begin{theorem}\label{thm:schur}
Assume that $m\in L^\infty(\RR^2)$ has a representative whose coordinate
sections are locally absolutely continuous away from $x=y$ for almost
every fixed value of the other coordinate, and let
\begin{equation}\label{eq:continuous-condition}
 M=\norm{m}_\infty,\qquad
 D_i=\esssup_{x\ne y}|x-y|\,|\partial_i m(x,y)|<\infty,
 \quad i=1,2.
\end{equation}
For every $X\in S_2(H)$ and $t>0$, the stronger operator inequalities hold:
\begin{align*}
 S_t((M_mX)^*M_mX)-(S_tM_mX)^*S_tM_mX
 &\le (M+2D_2)^2\bigl[S_t(X^*X)-(S_tX)^*S_tX\bigr],\\
 S_t((M_mX)(M_mX)^*)-(S_tM_mX)(S_tM_mX)^*
 &\le (M+2D_1)^2\bigl[S_t(XX^*)-(S_tX)(S_tX)^*\bigr].
\end{align*}
Consequently,
\begin{align}\label{eq:schur-result}
 \norm{M_mX}_{\bmo_c(S)}&\le(M+2D_2)\norm{X}_{\bmo_c(S)},\\
 \norm{M_mX}_{\bmo_r(S)}&\le(M+2D_1)\norm{X}_{\bmo_r(S)}.\nonumber
\end{align}
In particular, with $C_m=M+2\max(D_1,D_2)$,
\begin{equation}\label{eq:endpoint}
 \norm{M_mX}_{\bmo(S)}\le C_m\norm{X}_{\bmo(S)}
 \le C_m\norm{X}_\infty.
\end{equation}
All bounds hold at every matrix amplification.
\end{theorem}

\subsection{Proof by Schur tails}
\begin{lemma}[Continuous Schur-tail identity]\label{lem:schur-tails}
For fixed $t>0$, let $V_s^\sigma(X)$ have kernel
\[
 \mathbf1_{\{\sigma(x-y)\ge s\}}
 e^{-t(\sigma(x-y)-s)}K_X(x,y),\qquad s>0.
\]
For $X\in S_2(H)$,
\begin{equation}\label{eq:schur-variance}
 S_t(X^*X)-(S_tX)^*S_tX
 =2t\sum_{\sigma=\pm1}\int_0^\infty
 V_s^\sigma(X)^*V_s^\sigma(X)\,ds,
\end{equation}
with the positive integral converging in trace norm.
\end{lemma}
\begin{proof}
First let $k=K_X\in C_c^\infty(\RR^2)$, and write
\[
 W_t(X)=S_t(X^*X)-(S_tX)^*S_tX,\qquad
 v_s^\sigma(x,y)=\mathbf1_{\{\sigma(x-y)\ge s\}}
 e^{-t(\sigma(x-y)-s)}.
\]
By \eqref{eq:continuous-scalar},
\[
 e^{-t|y-z|}-e^{-t(|x-y|+|x-z|)}
 =2t\sum_{\sigma=\pm1}\int_0^\infty
 v_s^\sigma(x,y)v_s^\sigma(x,z)\,ds.
\]
Hence
\begin{align*}
 K_{W_t(X)}(y,z)
 &=\int_\RR\overline{k(x,y)}k(x,z)
 \bigl(e^{-t|y-z|}-e^{-t(|x-y|+|x-z|)}\bigr)\,dx\\
 &=2t\sum_\sigma\int_0^\infty\int_\RR
 \overline{v_s^\sigma(x,y)k(x,y)}
 v_s^\sigma(x,z)k(x,z)\,dx\,ds\\
 &=2t\sum_\sigma\int_0^\infty
 K_{V_s^\sigma(X)^*V_s^\sigma(X)}(y,z)\,ds.
\end{align*}
For arbitrary $X\in S_2(H)$, Tonelli gives
\begin{align*}
 2t\sum_\sigma\int_0^\infty\|V_s^\sigma(X)\|_2^2\,ds
 &=\int_{\RR^2}|k(x,y)|^2
 \left(2t\sum_\sigma\int_0^\infty|v_s^\sigma(x,y)|^2\,ds\right)dx\,dy\\
 &=\int_{\RR^2}(1-e^{-2t|x-y|})|k(x,y)|^2\,dx\,dy
 \le\|X\|_2^2.
\end{align*}
Thus the positive integral converges in $S_1$. Both sides of
\eqref{eq:schur-variance} are continuous from $S_2$ to $S_1$
by the above estimate and Schatten Cauchy--Schwarz; density completes
 the proof.
\end{proof}

\begin{proof}[Proof of Theorem~\ref{thm:schur}]
Fix $t>0$ and write $a_x^\sigma(u)=m(x,x-\sigma u)$. Let $A_s^\sigma$ and $B_u^\sigma$ be multiplication operators on $H$ with functions $a_x^\sigma(s)$ and $\partial_u a_x^\sigma(u)$, respectively. The hypotheses imply
\[
 \norm{A_s^\sigma}\le M,\qquad
 \norm{B_u^\sigma}\le D_2/u.
\]
The fundamental theorem of calculus along each coordinate section gives
\begin{equation}\label{eq:abel-schur}
 V_s^\sigma(M_mX)=A_s^\sigma V_s^\sigma(X)
 +\int_s^\infty e^{-t(u-s)}B_u^\sigma V_u^\sigma(X)\,du.
\end{equation}
To check the identity, fix $d=\sigma(x-y)\ge s$: the integral is over $s\le u\le d$ and its scalar coefficient telescopes from $m(x,x-\sigma s)$ to $m(x,y)$. For almost every $s>0$ the integral converges in $S_2$, by Cauchy--Schwarz and Lemma~\ref{lem:schur-tails}.

For each input vector $\eta\in H$, regard $(V_s^\sigma(X)\eta)_{\sigma,s}$ as a vector in
\[
 \mathcal K=L^2(\{\pm1\}\times\RR_+;H).
\]
For each fixed output coordinate $x$ and sign $\sigma$, \eqref{eq:abel-schur} acts on the $s$ variable as a multiplication operator of norm at most $M$, plus an integral operator dominated by $D_2\mathbf1_{\{u>s\}}/u$. The Hardy inequality \eqref{eq:hardy} bounds its norm by $M+2D_2$. Taking the direct integral over $x,\sigma$ gives the same bound on $\mathcal K$. Thus
\[
 \sum_\sigma\int_0^\infty V_s^\sigma(M_mX)^*V_s^\sigma(M_mX)\,ds
 \le(M+2D_2)^2
 \sum_\sigma\int_0^\infty V_s^\sigma(X)^*V_s^\sigma(X)\,ds.
\]
Lemma~\ref{lem:schur-tails} proves the column estimate. For $m^\dagger(x,y)=\overline{m(y,x)}$, one has $(M_mX)^*=M_{m^\dagger}(X^*)$ and $D_2(m^\dagger)=D_1(m)$, proving the row estimate. Scalar multiplication in $x$ commutes with matrix amplification, so every bound is completely bounded. Finally $\norm{X}_{\bmo(S)}\le\norm{X}_\infty$ follows from the variance formula and complete positivity.
\end{proof}

\paragraph{The $S_2$ framework.}
Both $X$ and $M_mX$ belong to $S_2(H)$, and the kernel products and tail
integrals above have their usual Hilbert--Schmidt and trace-class
meanings. Appendix~\ref{app:interpolation} derives the Schatten estimates
from precisely these $S_2$ bounds.

\paragraph{Difference symbols.}
If $m(x,y)=\mu(x-y)$ with
$\|\mu\|_\infty=M$ and $\esssup_{\xi\ne0}|\xi\mu'(\xi)|=D$,
then $D_1=D_2=D$. In particular, $m(x,y)=\operatorname{sgn}(x-y)$
has $M=1$ and $D_1=D_2=0$. Since $M_m^2X=X$ on $S_2(H)$, applying
the operator inequalities twice gives the familiar variance identity
\[
 S_t((M_mX)^*M_mX)-(S_tM_mX)^*S_tM_mX
 =S_t(X^*X)-(S_tX)^*S_tX.
\]
The row identity holds as well, so $M_m$ is an isometry for both Poisson
BMO norms. The proof above is direct on the matrix algebra and
does not use transference to Fourier multipliers.

\appendix
\section{Poisson semigroup BMO and quadratic interval BMO}\label{app:interval-bmo}

The equivalence proved below is well known. We include the proof for the
convenience of the reader and to specify the relation between the norm used
in this article and the quadratic interval $\BMO_2$ seminorm. All seminorms are homogeneous,
so functions that differ by a constant are identified.

For a locally square-integrable complex-valued function $f$ and a bounded interval
$I\subset\RR$, write $f_I=|I|^{-1}\int_I f$. Define
\[
 \|f\|_{\BMO_2}
 =\sup_I\left(\frac1{|I|}\int_I|f-f_I|^2\right)^{1/2}.
\]
Recall that $p_t(u)=t/(\pi(t^2+u^2))$.

\begin{theorem}\label{thm:interval-equivalence}
The spaces $\BMO_2(\RR)$ and $\bmo(P)$ coincide modulo constants, with
equivalent seminorms. More precisely,
\begin{equation}\label{eq:interval-poisson-bounds}
 \pi^{-1/2}\|f\|_{\BMO_2}
 \le\|f\|_{\bmo(P)}
 \le\sqrt{112/\pi}\,\|f\|_{\BMO_2}.
\end{equation}
Every $\BMO_2$ function satisfies
$\int_\RR |f(y)|^2(1+y^2)^{-1}\,dy<\infty$, so its Poisson averages
and variances are well defined. The numerical constants are not intended
to be optimal.
\end{theorem}

\subsection*{From Poisson variance to interval oscillation}
For any probability density $q$ for which the integrals exist, its mean
$a_q=\int qf$ satisfies the centering identity
\begin{equation}\label{eq:scalar-centering-appendix}
 \int q(y)|f(y)-a|^2\,dy
 =\int q(y)|f(y)-a_q|^2\,dy+|a_q-a|^2.
\end{equation}
In particular,
\[
 V_t^P(f)(x)
 =\int_\RR p_t(x-y)|f(y)-P_tf(x)|^2\,dy.
\]
Let $I=(x-t,x+t)$. On $I$, one has
$p_t(x-y)\ge(2\pi t)^{-1}=(\pi|I|)^{-1}$.
Since $f_I$ minimizes the unweighted square integral over constants,
\[
 \frac1{|I|}\int_I|f-f_I|^2
 \le\frac1{|I|}\int_I|f-P_tf(x)|^2
 \le\pi V_t^P(f)(x).
\]
Every bounded interval is of this form. Taking suprema proves the first
inequality in \eqref{eq:interval-poisson-bounds}.

\subsection*{From interval oscillation to Poisson variance}
Let $B=\|f\|_{\BMO_2}$, fix $x\in\RR$ and $t>0$, and put
\[
 I_j=(x-2^jt,x+2^jt),\quad a_j=f_{I_j},\qquad j\ge0.
\]
Cauchy--Schwarz on $I_{j-1}\subset I_j$ gives
\[
 |a_j-a_{j-1}|
 \le\left(\frac1{|I_{j-1}|}\int_{I_{j-1}}|f-a_j|^2\right)^{1/2}
 \le\sqrt2 B.
\]
Hence $|a_j-a_0|\le\sqrt2 jB$. The unweighted centering identity yields
\begin{equation}\label{eq:dyadic-mean-oscillation}
 \frac1{|I_j|}\int_{I_j}|f-a_0|^2
 =\frac1{|I_j|}\int_{I_j}|f-a_j|^2+|a_j-a_0|^2
 \le(1+2j^2)B^2.
\end{equation}
Set $A_0=I_0$ and $A_j=I_j\setminus I_{j-1}$ for $j\ge1$.
For $y\in A_j$ the Poisson density satisfies
\[
 p_t(x-y)\le\frac8\pi\,\frac{2^{-j}}{|I_j|},\qquad j\ge0.
\]
For $j=0$ this follows from $p_t\le1/(\pi t)$; for $j\ge1$ use
$|x-y|\ge2^{j-1}t$. It follows from
\eqref{eq:dyadic-mean-oscillation} that
\begin{align*}
 \int_\RR p_t(x-y)|f(y)-a_0|^2\,dy
 &\le\frac8\pi\sum_{j\ge0}2^{-j}
       \frac1{|I_j|}\int_{I_j}|f-a_0|^2\\
 &\le\frac8\pi B^2\sum_{j\ge0}2^{-j}(1+2j^2)
 =\frac{112}\pi B^2.
\end{align*}
This estimate first proves the required integrability: at $x=0,t=1$,
adding back the constant $a_0$ gives
$\int |f(y)|^2(1+y^2)^{-1}\,dy<\infty$.
For any other fixed $x,t$, the weight $p_t(x-y)$ is comparable to
$(1+y^2)^{-1}$, so $P_tf(x)$ and $P_t|f|^2(x)$ exist as well.
Now \eqref{eq:scalar-centering-appendix}, with $a=a_0$, bounds
$V_t^P(f)(x)$ by the preceding integral. Taking suprema proves the second
inequality in \eqref{eq:interval-poisson-bounds}, and completes the proof
of Theorem~\ref{thm:interval-equivalence}.

\begin{remark}[Matrix-valued quadratic interval BMO]
For a matrix-valued function $F$, define
\[
 \|F\|_{\BMO_{2,c}^{\mathrm{int}}}
 =\sup_I\left\|\frac1{|I|}\int_I(F-F_I)^*(F-F_I)\right\|^{1/2}.
\]
The same comparison holds between this norm and $\bmo_c(P)$, with the
constants in \eqref{eq:interval-poisson-bounds} independent of matrix size.
Indeed, the two square-oscillation arguments apply to the Hilbert-space-valued
functions $F\eta$, uniformly for unit vectors $\eta$. Apply them to $F^*$
for the row norm, and take the maximum for full matrix-valued BMO.
\end{remark}

\section{Interpolation}\label{app:interpolation}

We pass directly from the BMO and Hilbert--Schmidt estimates to Schatten
classes. We use the semigroup BMO framework of Mei \cite{MeiTent} and the
semigroup-to-martingale interpolation theorem of Junge--Mei \cite{JM}.
The endpoint proofs in Sections~\ref{sec:fourier} and~\ref{sec:schur} are
independent of that theorem. We spell out its application to the present
Schur semigroup and the resulting multiplier estimate.

\subsection{The interpolation identity}

Let $\mathcal B=\bmo(S)$ be the full row-and-column Poisson space in the
compatible completed Hilbert-module realization of
\cite[Section 5.2.2]{JM}. On bounded operators its norm is
\[
 \|X\|_{\mathcal B}
 =\max\left\{\sup_{t>0}\|S_t(X^*X)-(S_tX)^*S_tX\|^{1/2},
 \sup_{t>0}\|S_t(XX^*)-(S_tX)(S_tX)^*\|^{1/2}\right\}.
\]
Multiplication operators form the null space. The completed column and row
spaces are identified through their common Poisson derivatives
$\partial_t S_tX$, as in the cited construction. This specifies a compatible
couple with $S_2(H)$; it is larger than the BMO norm completion of $S_2(H)$.
The interpolation brackets below denote the complex method.

\begin{theorem}\label{thm:schur-interpolation}
For $2<p<\infty$,
\begin{equation}\label{eq:bmo-s2-interpolation}
 [\mathcal B,S_2(H)]_{2/p}=S_p(H).
\end{equation}
There is an absolute constant $C$ such that
\begin{equation}\label{eq:interp-constants}
 \|X\|_{[\mathcal B,S_2]_{2/p}}
 \le\|X\|_{S_p}
 \le Cp\,\|X\|_{[\mathcal B,S_2]_{2/p}}.
\end{equation}
The constants are uniform in matrix amplification.
\end{theorem}

\begin{proof}
The underlying heat Schur semigroup is
\[
 K_{H_sX}(x,y)=e^{-s|x-y|^2}K_X(x,y),
\]
whose subordinate Poisson semigroup is $S$. If $(B_s)$ is standard real
Brownian motion and $(U_v\eta)(x)=e^{ivx}\eta(x)$, then
\[
 \pi_s(X)(\omega)=U_{\sqrt2 B_s(\omega)}XU_{\sqrt2 B_s(\omega)}^*
\]
is a trace-preserving $*$-homomorphism into
$L^\infty(\Omega)\,\overline\otimes\,B(H)$. Independence of increments gives
\[
 \mathbb E_s\pi_u(X)=\pi_s(H_{u-s}X),\qquad 0\le s\le u.
\]
Thus the underlying heat semigroup admits a standard Markov dilation.

We use the proof of \cite[Theorem 5.2(i)]{JM} with base exponent $2$,
in the semifinite realization of \cite[Section 5.2.2 and Theorem 5.12]{JM}.
The argument passes through Meyer's stopped-martingale model: its
semigroup-to-martingale maps are uniformly bounded on the Poisson BMO
endpoint and on $L_2$, and the recovery estimate returns the $L_p$ norm
of the original element from the martingale model. Applying the
martingale identity
\[
 [\mathrm{BMO}_{\mathrm{mart}},L_2]_{2/p}=L_p
\]
inside that model yields
\[
 \|X\|_{S_p^0}
 \le Cp\,\|X\|_{[\mathcal B,S_2^0]_{2/p}}.
\]
Here the superscript $0$ removes the fixed part. This is precisely the
direct BMO--$L_2$ argument of the cited theorem; no trace-class endpoint
is used in this application.

In our setting the reduced finite Schatten spaces are the full ones.
Indeed, on $S_2$, dominated convergence applied to the kernel gives
$H_sX\to0$ in $S_2$ as $s\to\infty$, since the diagonal has measure zero.
For $p>2$ and $X\in S_2$, contractivity on $B(H)$ and
\[
 \|H_sX\|_{S_p}
 \le\|H_sX\|_{S_2}^{2/p}\|H_sX\|^{1-2/p}
\]
give convergence in $S_p$. Finite-rank density and contractivity then give
the same convergence for every $X\in S_p$. Thus $S_2^0=S_2$ and $S_p^0=S_p$.
The same conclusions hold for the subordinate Poisson semigroup.

For the other norm inequality, interpolate the contractive map
$B(H)\to\mathcal B$ with the identity on $S_2$ and use
$[B(H),S_2]_{2/p}=S_p$.
All constructions commute with matrix amplification.

\end{proof}

\begin{corollary}[Interpolation from the $S_2$ domain]\label{cor:operator-interpolation}
Suppose a linear map $T:S_2(H)\to S_2(H)$ satisfies
\[
 \|TX\|_{\bmo(S)}\le A\|X\|_{\bmo(S)},\qquad
 \|TX\|_{S_2}\le B\|X\|_{S_2},\qquad X\in S_2(H).
\]
Then $T$ extends to $S_p(H)$ and
\[
 \|TX\|_{S_p}\le Cp\,A^{1-2/p}B^{2/p}\|X\|_{S_p},
 \qquad 2<p<\infty.
\]
The assertion is uniform in matrix amplification whenever the endpoint
bounds are uniform.
\end{corollary}
\begin{proof}
Let $\mathcal B_0$ be the closure of $S_2(H)$ in the compatible space
$\mathcal B$ above. Since
$\|X\|_{\mathcal B}\le\|X\|\le\|X\|_{S_2}$, we have
$\mathcal B\cap S_2=S_2$. The regularization theorem for complex
interpolation \cite[Section 9.3]{Calderon} therefore gives, isometrically,
\[
 [\mathcal B_0,S_2]_\theta=[\mathcal B,S_2]_\theta,
 \qquad 0<\theta<1.
\]
Here regularization means replacing each endpoint by the closure of
its intersection with the other endpoint. Thus the BMO norm completion
of $S_2$ is sufficient for interpolation.

The first hypothesis extends $T$ by continuity to $\mathcal B_0$ with
norm at most $A$, compatibly with its $S_2$ action. Put $\theta=2/p$.
Exactness of complex interpolation bounds $T$ on
$[\mathcal B_0,S_2]_\theta$ by $A^{1-\theta}B^\theta$.
Apply \eqref{eq:interp-constants} and the regularization identity to
obtain the assertion. The extension to $S_p$ is unique by finite-rank
density.
\end{proof}

\subsection{Schatten estimates for Schur multipliers}
Theorem~\ref{thm:schur} proves on $S_2(H)$ that
\[
 \|M_mX\|_{\bmo(S)}\le A\|X\|_{\bmo(S)},\qquad
 A=M+2\max\{D_1,D_2\}.
\]
On $S_2$, multiplication of kernels gives $\|M_m\|_{S_2\to S_2}\le M$.
Corollary~\ref{cor:operator-interpolation} therefore gives
\[
 \|M_mX\|_{S_p}\le Cp\,A^{1-2/p}M^{2/p}\|X\|_{S_p},\qquad p>2.
\]
For $1<p<2$, use trace duality and the symbol $\overline m$, whose
derivative bounds are unchanged. With $p^*=\max\{p,p/(p-1)\}$ we obtain
\[
 \|M_mX\|_{S_p}
 \le Cp^*\bigl(M+2\max\{D_1,D_2\}\bigr)\|X\|_{S_p},
 \qquad 1<p<\infty.
\]
At $p=2$ the direct bound is $M$. Every estimate remains valid at all
matrix amplifications.

\subsection{The Fourier multiplier consequence}
The parallel interpolation statement on
$L^\infty(\RR)\,\overline\otimes\,\cM$, for a semifinite coefficient algebra
$(\cM,\tau)$, is
\[
 [\bmo(P),L^2(\RR;\cM)]_{2/p}=L^p(\RR;\cM),\qquad p>2,
\]
with the full row-and-column endpoint in the same compatible realization.
The underlying heat semigroup has the Brownian dilation
$\pi_s f(x,\omega)=f(x+\sqrt2B_s(\omega))$.
The preceding martingale argument therefore applies. Its reduced finite
$L^p$ spaces are the full spaces: the heat semigroup tends to zero on
$L^2$ by Plancherel, and contractivity and density give the same conclusion
for $2<p<\infty$.

To apply this identity using only the $L^2$ domain, put
\[
 E=L^2(\RR;\cM)\cap\bmo(P),\qquad
 \mathcal B_{P,0}=\overline E^{\,\bmo(P)}.
\]
The intersection $E$ is dense in $L^2$: finite sums of scalar Schwartz
functions with coefficients in $\cM\cap L^2(\cM,\tau)$ already form a
dense subspace. Regularization as above gives
\[
 [\mathcal B_{P,0},L^2]_{2/p}=[\bmo(P),L^2]_{2/p}=L^p.
\]
Corollary~\ref{cor:matrix-fourier} shows that $T_\mu$ preserves $E$ and
extends to $\mathcal B_{P,0}$ with norm at most $M+2D$.
Together with its Plancherel bound $M$, complex interpolation gives
\[
 \|T_\mu f\|_p
 \le Cp\,(M+2D)^{1-2/p}M^{2/p}\|f\|_p,\qquad p>2.
\]
Duality with the symbol $\overline\mu$ then yields
\[
 \|T_\mu f\|_p\le Cp^*(M+2D)\|f\|_p,
 \qquad 1<p<\infty.
\]
For $\cM=M_q$, these are the usual matrix-valued $L^p$ norms, with constants
independent of $q$. The scalar case is recovered by taking $\cM=\mathbb C$.

\section{A direct Hardy-space estimate}\label{app:hardy}
Let $H_P^1(\RR)$ denote the usual real Hardy space associated with the
Poisson semigroup. We use its equivalent square-function norm
\[
 G(f)(x)=\left(\int_0^\infty t\,|\partial_tP_tf(x)|^2\,dt\right)^{1/2},
 \qquad \|f\|_{H_P^1}=\|G(f)\|_1.
\]
The square-function characterization and $H^1$--BMO duality are classical;
see \cite{FS}. We also use the standard $L^2$ atomic characterization
\cite{Coifman}: an atom is supported in an interval $I$, has integral zero,
and has $L^2$ norm at most $|I|^{-1/2}$. The infimum of the sums of the
absolute values of the coefficients in atomic representations is
equivalent to the displayed Hardy norm.

\begin{theorem}\label{thm:hardy-fourier}
Under the hypotheses of Theorem~\ref{thm:real-fourier}, the multiplier
$T_\mu$ extends to a bounded operator on $H_P^1(\RR)$, with
\[
 \|T_\mu f\|_{H_P^1}\le C(M+D)\|f\|_{H_P^1},
\]
where $C$ is absolute. More generally, if
$a,(x-c)a\in L^2(\RR)$ and $\int_\RR a=0$, then
\begin{equation}\label{eq:hardy-molecule}
 \|T_\mu a\|_{H_P^1}
 \le C(M+D)\bigl(\|a\|_2\,\|(x-c)a\|_2\bigr)^{1/2}.
\end{equation}
No support assumption is needed in this latter estimate.
\end{theorem}

\begin{proof}
Translate so that $c=0$, and set $A=\|a\|_2$, $B=\|xa\|_2$.
Put $u=\widehat a$, using the Fourier convention of
Section~\ref{sec:fourier}. Weighted Cauchy--Schwarz gives
$a\in L^1$, and cancellation gives $u(0)=0$. Hardy's inequality on the
positive and negative half-lines, obtained by integration by parts,
therefore gives
\begin{equation}\label{eq:frequency-hardy}
 \|u/\xi\|_2\le2\|u'\|_2=4\pi B.
\end{equation}
With $\rho=2\pi|\xi|$, Plancherel and $\int_0^\infty t\rho^2e^{-2t\rho}\,dt=1/4$ give
\begin{equation}\label{eq:hardy-unweighted}
 \|G(T_\mu a)\|_2\le \tfrac12 M A.
\end{equation}

For the weighted estimate, multiplication by $x$ corresponds to
$(i/(2\pi))\partial_\xi$. Differentiation of
$\rho e^{-t\rho}\mu(\xi)u(\xi)$ gives
\[
 2\pi\operatorname{sgn}(\xi)(1-t\rho)e^{-t\rho}\mu u
 +\rho e^{-t\rho}(\mu'u+\mu u').
\]
There is no boundary term at zero: $u(0)=0$ and $\mu$ is bounded.
Using
\[
 \int_0^\infty t(1-t\rho)^2e^{-2t\rho}\,dt=\frac1{8\rho^2},
 \qquad
 \int_0^\infty t\rho^2e^{-2t\rho}\,dt=\frac14,
\]
Minkowski and Plancherel yield
\begin{align}
 \|xG(T_\mu a)\|_2
 &\le\left(\frac{M}{4\pi\sqrt2}+\frac D{4\pi}\right)\|u/\xi\|_2
       +\frac M{4\pi}\|u'\|_2 \notag\\
 &\le\left[\left(\frac1{\sqrt2}+\frac12\right)M+D\right]B.
 \label{eq:hardy-weighted}
\end{align}
All differentiations are valid weakly; \eqref{eq:frequency-hardy} controls
the possibly singular term $\mu'u$.

For every $r>0$, weighted Cauchy--Schwarz now gives
\begin{align*}
 \|G(T_\mu a)\|_1
 &\le(\pi r)^{1/2}
 \left(\|G(T_\mu a)\|_2^2+
             r^{-2}\|xG(T_\mu a)\|_2^2\right)^{1/2}\\
 &\le C(M+D)(rA^2+r^{-1}B^2)^{1/2}.
\end{align*}
For $a\ne0$, choose $r=B/A$ to prove \eqref{eq:hardy-molecule};
the zero case is immediate.

If $a$ is an atom supported in $(c-r,c+r)$, then
$\|(x-c)a\|_2\le r\|a\|_2$ and
$\|a\|_2\le(2r)^{-1/2}$. Hence \eqref{eq:hardy-molecule} bounds
$\|T_\mu a\|_{H_P^1}$ by $C(M+D)$, uniformly over atoms.
Summing an atomic representation proves the theorem.
The extension agrees with the Fourier multiplier in distributions:
atomic sums converge in $H_P^1$, hence in $L^1$, so their Fourier
transforms converge uniformly, and multiplication by the bounded
function $\mu$ preserves the resulting distributional limit.
\end{proof}

The calculation also applies to $\partial_xP_t$, whose Fourier factor is
$2\pi i\xi e^{-2\pi t|\xi|}$. It thus proves the same estimates for the equivalent
full-gradient square function. In particular, mean-zero functions with
$\|a\|_2\le r^{-1/2}$ and $\|(x-c)a\|_2\le r^{1/2}$ satisfy the same
uniform bound even when they have noncompact support.

\begin{corollary}[Operator-valued Hardy spaces]\label{cor:hardy-matrix}
Let $(\cM,\tau)$ be a semifinite von Neumann algebra. Define
\[
 G_c^\nabla(f)(x)
 =\left(\int_0^\infty t\left[
   (\partial_tP_tf)^*(\partial_tP_tf)
  +(\partial_xP_tf)^*(\partial_xP_tf)\right](x)\,dt\right)^{1/2},
\]
and use the column norm
$\|f\|_{H_{P,c}^1}=\int_\RR\tau(G_c^\nabla(f)(x))\,dx$.
Then
\[
 \|T_\mu f\|_{H_{P,c}^1}\le C(M+D)\|f\|_{H_{P,c}^1},
\]
with an absolute constant independent of $\cM$. The corresponding
statement holds for the row space and for the sum of the column and row
Hardy spaces.
\end{corollary}
\begin{proof}
Use the standard square-function and column atomic characterizations
\cite[Corollary 2.7 and Theorem 2.8]{MeiHardy}.
For a column atom supported in $(c-r,c+r)$, put
$A=(\int_\RR a^*a)^{1/2}$, so that $\tau(A)\le(2r)^{-1/2}$.
First take $A$ bounded; general atoms follow by right spectral truncation.
Apply the Hilbert-space-valued versions of
\eqref{eq:hardy-unweighted} and \eqref{eq:hardy-weighted} to $a(\cdot)\eta$
for each $\eta\in\HH$, for both derivatives. This gives, in operator order,
\[
 \int_\RR \left(1+\frac{(x-c)^2}{r^2}\right)
      G_c^\nabla(T_\mu a)(x)^2\,dx
 \le C(M+D)^2 A^2.
\]
Operator concavity of the square root, with the probability density
$(\pi r)^{-1}(1+(x-c)^2/r^2)^{-1}$, gives
\[
 \int_\RR G_c^\nabla(T_\mu a)(x)\,dx
 \le\sqrt{\pi r}
 \left(\int_\RR\left(1+\frac{(x-c)^2}{r^2}\right)
             G_c^\nabla(T_\mu a)(x)^2\,dx\right)^{1/2}.
\]
Taking the trace proves the uniform atom bound.
Sum the atoms, and use adjoints for the row estimate.
\end{proof}


\begin{thebibliography}{8}
\bibitem{CGPT}
J. M. Conde-Alonso, A. M. Gonz\'alez-P\'erez, J. Parcet and E. Tablate,
\emph{Schur multipliers in Schatten--von Neumann classes},
Annals of Mathematics (2) \textbf{198} (2023), no. 3, 1229--1260.
\href{https://doi.org/10.4007/annals.2023.198.3.5}{doi:10.4007/annals.2023.198.3.5}.
\bibitem{MeiTent}
T. Mei,
\emph{Tent spaces associated with semigroups of operators},
Journal of Functional Analysis \textbf{255} (2008), no. 12, 3356--3406.
\href{https://doi.org/10.1016/j.jfa.2008.09.021}{doi:10.1016/j.jfa.2008.09.021}.
\bibitem{JM}
M. Junge and T. Mei,
\emph{BMO spaces associated with semigroups of operators},
Mathematische Annalen \textbf{352} (2012), no. 3, 691--743.
\href{https://doi.org/10.1007/s00208-011-0657-0}{doi:10.1007/s00208-011-0657-0}.
Also available as \href{https://arxiv.org/abs/1111.6194}{arXiv:1111.6194}.
We use the section and theorem numbering of that preprint.
\bibitem{CMR}
I. Cano, T. Mei and \'E. Ricard, \emph{Schur multipliers on $H_1^c$}, forthcoming preprint.
\bibitem{FS}
C. Fefferman and E. M. Stein,
\emph{$H^p$ spaces of several variables},
Acta Mathematica \textbf{129} (1972), 137--193.
\href{https://doi.org/10.1007/BF02392215}{doi:10.1007/BF02392215}.
\bibitem{Coifman}
R. R. Coifman,
\emph{A real variable characterization of $H^p$},
Studia Mathematica \textbf{51} (1974), no. 3, 269--274.
\href{https://doi.org/10.4064/sm-51-3-269-274}{doi:10.4064/sm-51-3-269-274}.
\bibitem{MeiHardy}
T. Mei, \emph{Operator valued Hardy spaces},
Memoirs of the American Mathematical Society \textbf{188} (2007), no. 881.
Also available as \href{https://arxiv.org/abs/math/0311160}{arXiv:math/0311160}.
The theorem numbering used here refers to that preprint.

\bibitem{Calderon}
A. P. Calder\'on,
\emph{Intermediate spaces and interpolation, the complex method},
Studia Mathematica \textbf{24} (1964), 113--190.
\href{https://doi.org/10.4064/sm-24-2-113-190}{doi:10.4064/sm-24-2-113-190}.

\end{thebibliography}
\end{document}